\documentclass[letterpaper, 10 pt, conference]{ieeeconf}  

\IEEEoverridecommandlockouts                              
\usepackage{epsfig} 
\usepackage{amsmath,amssymb,amsfonts}
\usepackage{color}
\usepackage{booktabs}
\usepackage{cite}
\usepackage{multirow}
\usepackage{mathtools}
\usepackage{pdfrender}

\usepackage{centernot}
\usepackage{color}
\usepackage{subfig}
\usepackage[acronym]{glossaries}

\makeglossaries
\newacronym{LP}{LP}{Linear Programming}
\newacronym{MPC}{MPC}{Model Predictive Control}
\newacronym{QP}{QP}{Quadratic Programming}
\newacronym{LTI}{LTI}{Linear Time Invariant}
\newacronym{LQR}{LQR}{Linear Quadratic Regulator}
\newacronym{ARE}{ARE}{Algebraic Riccati Equation}
\newacronym{CLF}{CLF}{Control Lyapunov Function}
\newacronym{QCLF}{QCLF}{Quadratic Control Lyapunov Function}
\newacronym{PQCLF}{PQCLF}{Piecewise Quadratic Control Lyapunov Function}
\newacronym{CSCLF}{CSLF}{Constraints-Shaped Lyapunov Function}
\newacronym{MISO}{MISO}{Multiple Input - Single Output}
\newacronym{SISO}{SISO}{Single Input - Single Output}
\newacronym{MIMO}{MIMO}{Multiple Input - Multiple Output}
\newacronym{EAS}{EAS}{Euler Auxiliary System}
\newacronym{DA}{DA}{Domain of Attraction}
\newacronym{psd}{p.s.d.}{positive semi-definite}
\makeglossaries
\newacronym{MIQP}{MIQP}{Mixed-Integer Quadratic Programming}
\newacronym{BSI}{BSI}{Bilateral Symmetric Interaction}
\newacronym{AD}{AD}{Automated Driving}
\newacronym{MILP}{MILP}{Mixed-Integer Linear Programming}
\newacronym{MLD}{MLD}{Mixed-Logical-Dynamical}
\newacronym{GS}{GS}{Gauss-Southwell}
\newacronym{GNEP}{GNEP}{Generalized Nash Equilibrium Problem}
\newacronym{NEP}{NEP}{Nash Equilibrium Problem}
\newacronym{NE}{NE}{Nash Equilibrium}
\newacronym{MINE}{$\varepsilon$-MINE}{$\varepsilon$-Mixed-Integer Nash Equilibrium}

\graphicspath{{Pictures/}}

\newcommand{\R}{\mathbb{R}}

\newcommand{\bbS}{\mathbb{S}}
\newcommand{\mc}[1]{\mathcal{#1}}

\newcommand{\diag}{\mathrm{diag}}

\newcommand{\mb}[1]{\mathbf{#1}}
\newcommand{\bs}[1]{\boldsymbol{#1}}

\newcommand{\continuanceref}{}

\newtheorem{problem}{Problem}
\newtheorem{definition}{Definition}
\newtheorem{proposition}{Proposition}
\newtheorem{lemma}{Lemma}
\newtheorem{corollary}{Corollary}

\newtheorem{standing}{Standing Assumption}

\title{\LARGE \bf
Nash equilibrium seeking in potential games\\ with double-integrator agents
}

\author{Filippo Fabiani and
	 Andrea Caiti
\thanks{The authors are with the Department of Information Engineering, University of Pisa, Italy
	{\tt\small (filippo.fabiani@unipi.it, andrea.caiti@unipi.it)}.}%
}

\begin{document}

\maketitle
\thispagestyle{empty}
\pagestyle{empty}

\begin{abstract}
In this paper, we show the equivalence between a constrained, multi-agent control problem, modeled within the port-Hamiltonian framework, and an exact potential game. Specifically, critical distance-based constraints determine a network of double-integrator agents, which can be represented as a graph. Virtual couplings, i.e., pairs of spring-damper, assigned to each edge of the graph, allow to synthesize a distributed, gradient-based control law that steers the network to an invariant set of stable configurations. We characterize the points belonging to such set as Nash equilibria of the associated potential game, relating the parameters of the virtual couplings with the equilibrium seeking problem, since they are crucial to shape the transient behaviour (i.e., the convergence) and, ideally, the set of reachable equilibria.
\end{abstract}

\section{Introduction}
Distributed control of networked, multi-agent systems is a hot topic within the system-and-control community, since a peculiar characteristic of modern society is the ubiquitousness of large-scale systems with a complex network structure, involving interacting, (possibly) autonomous subsystems.

In the last few years, the control of agents governed by a double-integrator
dynamics has been widely investigated, both including typical
consensus protocols \cite{cheng2011necessary,li2011finite,qin2012sufficient,dong2013leader} and formation/distance-based issues \cite{hao2013stability,oh2014distance,fabiani2016distributed,fabiani2016passivity}. 
In this context, we consider double-integrator agents which have to satisfy critical, distance-based constraints that couple each others, defining a network that can be described as vertices and edges of an arbitrarily oriented graph. Successively, the network is recast within the port-Hamiltonian framework, closely related with the passivity theory \cite{duindam2009modeling,van2014port}. Virtual couplings, i.e., pairs of spring and damper in parallel, are assigned to each edge, in order to define the interaction forces among neighboring agents, while preserving the passivity of the network. This leads to a distributed, gradient-based control law that steers the system to a stable configuration. 

On the same leitmotiv of \cite{fabiani2018distributed}, we identify an intriguing affinity between the port-Hamiltonian formulation of the constrained control problem and a potential game, which revolves around the concept of \gls{NE}. 

Along this direction, the convergence to a \gls{NE} in games involving continuous-time, passive systems has been recently addressed in several works. Specifically, \cite{gadjov2017continuous,gadjov2018passivity,gao2018passivity} proposed passive techniques to solve the Nash equilibrium seeking problem over networks and in finite games, while \cite{fox2013population,mabrok2016passivity,park2018passivity} focused on the relation between passivity and evolutionary/stable games. 

Our work moves towards a novel interpretation of physical, multi-agent systems admitting a port-Hamiltonian model, providing a potential game-theoretic perspective. We envision that each agent aims at minimizing its energetic contribution within the network by seeking for an ad hoc strategy that satisfies distance-based constraints.
Here, the role of the virtual couplings is key, since they exhibit symmetries across the decision variables of the agents, i.e., two connected agents showing the same deviation in term of strategy reflect in exactly the same amount of deviation on the respective objective functions. This is crucial to prove the existence of an exact potential game associated to the related control problem. Moreover, the parameters of each pair of spring-damper can be chosen to shape both the transient behaviour of the network and the set of reachable equilibria within the associated \gls{NEP} \cite{facchinei2007generalized}.

The paper is organized as follows: after some basic preliminary recall (\S\ref{sec:preliminaries}), we model the constrained, multi-agent control problem within the port-Hamiltonian framework, synthesizing a distributed, gradient-based control law that steers the network to a stable configuration (\S\ref{sec:controlprob_pH}). Successively, we attach a potential game-theoretic perspective to the addressed control problem, providing equivalence results that characterize the set of Nash equilibria and the convergence to one of them (\S\ref{sec:potgamechar}). Finally, numerical simulations support the theoretical results given in the previous sections (\S\ref{sec:num_sim}).
\subsubsection*{Notation} 
$\R$, $\R_{> 0}$ and $\R_{\geq 0}$ denote the set of real, positive real, non-negative real numbers, respectively. $\bbS^{n}_{\succ 0}$ ($\bbS^{n}_{\succcurlyeq 0}$) denotes the set of symmetric, positive (semi-)definite matrices. Given vectors $x_1, \ldots, x_N \in \R^n$, $\bs{x} \coloneqq \left(x_1; \ldots; x_N\right)$ denotes $\left(x^\top_1,\ldots,x^\top_N\right)^\top \in \R^{nN}$. $A \otimes B$ denotes the Kronecker product between matrices $A$ and $B$. $\|x\|$ is the 2-norm of the vector $x$, while $|\mc{S}|$ denotes the cardinality of the set $\mc{S}$. $\mc{C}^k$ is the class of $k$-times continuously differentiable functions.

\section{Preliminaries}\label{sec:preliminaries}
We start with some basic definitions of port-Hamiltonian systems and game theory.
Specifically, the generalized input-state-output dynamics within the port-Hamiltonian framework reads as \cite{duindam2009modeling}:
\begin{equation*}
\left\{\begin{aligned}
&\dot{x} = (S(x) - D(x)) \frac{\partial{H}}{\partial{x}} + G(x) u,\\
&y = G^{\top}(x) \frac{\partial{H}}{\partial{x}},
\end{aligned}
\right.
\end{equation*}
where $x \in \R^{n}$ denotes the state, $u \in \R^{m}$ the control input and $y \in \R^{m}$ the output. $S \in \R^{n \times n}$, $D \in \bbS^{n}_{\succcurlyeq 0}$ and $G \in \R^{n \times m}$ are a skew-symmetric, a dissipation and an input matrix, respectively, while $H:\R^{n} \to \R_{\geq 0}$ is the Hamiltonian function. A system admitting a port-Hamiltonian representation is passive with storage function $H$, since it directly falls into the following definition. 

\smallskip
\begin{definition}[\hspace{-0.01em}\cite{duindam2009modeling}]\label{def:passivity}
	A map $u \mapsto y$ is passive if there exists a $\mc{C}^1$, lower bounded function of the state, $V:\R^{n} \to \R_{\geq 0}$ (storage function), such that
	\begin{equation*}
		\dot{V}(x) \leq u^{\top} y \iff V(x(t)) - V(x(0)) \leq \int_{0}^{t} u^{\top}(\tau) y(\tau) \, d\tau.
	\end{equation*}
	\vspace{-.2cm}{\hfill$\square$}
\end{definition}
\smallskip

In this paper, we show some equivalence results that connect a constrained, multi-agent control problem, modeled within the port-Hamiltonian framework, and a potential game. With this aim, a game 
$\Gamma \coloneqq (\mc{I}, \{J_i\}_{i \in \mc{I}}, \{\mc{X}_i\}_{i \in \mc{I}})$
consists of $N$ agents, indexed by the set $\mc{I} \coloneqq \{1, \ldots, N\}$, each one controlling its own variable, $x_i \in \mc{X}_i \subseteq \R^{n}$, and aims at minimizing its objective function, $J_i:\R^{nN} \to \R$. Hence, we refer to $\bs{x} \coloneqq \left(x_1; \ldots; x_N\right) \in \R^{nN}$ as the collective vector of strategies and to $\bs{x}_{-i} \in \R^{(n - 1)N}$ as the vector of all the players' decisions except those of player $i$. To emphasize the $i$-th decision variable within the collective vector, sometimes we write $\bs{x}$ as $\left(x_i, \bs{x}_{-i}\right)$.

%
For the remainder of this section, we assume that, for all $i \in \mc{I}$ and all $\bs{x}_{-i}$, $J_i \in \mc{C}^1$, $J_i(\cdot, \bs{x}_{-i})$ is convex and $\mc{X}_i$ is closed and convex. Thus, potential games \cite{monderer1996potential} belong to a particular class of games characterized by the existence of a potential function $P : \bs{\mc{X}} \to \R$, with $\bs{\mc{X}} \coloneqq \prod_{i \in \mc{I}} \mc{X}_i$, such that, for all $i \in \mc{I}$, for all $\bs{x}_{-i}$, and for all $x_i, y_i \in \mc{X}_i$
\begin{equation}\label{eq:exactpotgame}
J_i(x_i, \bs{x}_{-i}) - J_i(y_i, \bs{x}_{-i}) = P(x_i, \bs{x}_{-i}) - P(y_i, \bs{x}_{-i}).
\end{equation}
One of the key ingredients of game theory is the concept of \gls{NE}, defined as follows.

\smallskip
\begin{definition}
	A collective vector $\bs{x}^{\ast} \coloneqq \left(x_1^\ast; \ldots; x_N^\ast\right) \in \bs{\mc{X}}$ is a Nash equilibrium of the game $\Gamma$ if, for all $i \in \mc{I}$,
	\begin{equation}\label{eq:NEP}
	J_i(x_i^{\ast}, \bs{x}^{\ast}_{-i}) \leq \underset{x_i \in \mc{X}_i}{\textrm{inf}} J_i(x_i, \bs{x}^{\ast}_{-i}).
	\end{equation}
	\vspace{-.2cm}{\hfill$\square$}
\end{definition}
\smallskip

Since each feasible set $\mc{X}_i$ is independent from the strategies adopted by the neighbors, $\bs{x}_{-i}$, the problem in \eqref{eq:NEP} generically refers to a \gls{NEP} \cite{facchinei2007generalized}.
We finally introduce the pseudo-gradient mapping of the game, $\mb{F} : \bs{\mc{X}} \to \R^n$, defined as $\mb{F}(\bs{x}) \coloneqq \left(\nabla_{x_i} J_i(\bs{x})\right)_{i \in \mc{I}}$, and we characterize a \gls{NE} via variational inequalities, according to the following result.

\smallskip
\begin{lemma}[\hspace{-0.01em}{\cite[Cor.~1]{facchinei2007generalized}}]
	A collective vector $\bs{x}^{\ast} \in \bs{\mc{X}}$ is a variational equilibrium of the game $\Gamma$ if and only if it satisfies the variational inequality
	\begin{equation*}
	(\bs{y} - \bs{x}^{\ast})^\top \mb{F}(\bs{x}^{\ast}) \geq 0, \quad \text{for all } \bs{y} \in \bs{\mc{X}}.
	\end{equation*}
	\vspace{-.2cm}{\hfill$\square$}
\end{lemma}
\section{Constrained control of agents with double-integrator dynamics}\label{sec:controlprob_pH}
This section introduces the constrained, multi-agent control problem addressed and derives the port-Hamiltonian model for the set of double-integrator agents. Successively, virtual couplings are introduced to control the overall system, forcing the distance-based constraints that couple the agents.
In details, we deal with the following problem.
\smallskip
\begin{problem}
	A set of double-integrator agents has to be controlled to a stable configuration, while satisfying some critical, distance-based constraints. 
	{\hfill$\square$}
\end{problem}
\smallskip

Specifically, we consider $N$ agents belonging to the set $\mc{I} \coloneqq \{1, \ldots, N\}$, where each agent $i \in \mc{I}$ is a single point (unitary) mass governed by a double-integrator dynamics:
\begin{equation*}
\ddot{q}_i = u_i.
\end{equation*}

Here, $q_i \in \R^n$ denotes the generalized coordinate of the $i$-th agent, and $u_i \in \R^n$ its control input. The corresponding linear momentum $p_i \in \R^{n}$, is given by
$
p_i = \dot{q}_i,
$
for all $i \in \mc{I}$. By recalling that the kinetic energy associated to the $i$-th mass, $h^{\textrm{k}}_i:\R^n \to \R_{\geq 0}$, is explicitly given by 
$$h^{\textrm{k}}_i(\dot{q}_i) = \frac{1}{2} \dot{q}_i^\top \dot{q}_i = \frac{1}{2} p_i^\top p_i,$$ 
and by assuming that each output $y_i = \tfrac{\partial h^{\textrm{k}}_i}{\partial p_i}$, $i \in \mc{I}$, the single agent dynamics in the port-Hamiltonian framework reads as:
\begin{equation*}
\forall i \in \mc{I} \; : \; 
\left(
\begin{array}{c}
\dot{q}_i\\
\dot{p}_i
\end{array}
\right) = 
\left(
\begin{array}{cc}
0 & I_n\\
0 & 0
\end{array}\right)
\left(
\begin{array}{c}
\frac{\partial h^{\textrm{k}}_i}{\partial q_i}\\
\frac{\partial h^{\textrm{k}}_i}{\partial p_i}
\end{array}
\right) + 
\left(
\begin{array}{c}
0\\
I_n
\end{array}
\right) u_i.
\end{equation*}

Then, to compactly characterize the dynamics of the whole set, we introduce $q \coloneqq (q_1; \ldots; q_N) \in \R^{n N}$ as the vector of generalized coordinates, and $u \coloneqq (u_1; \ldots; u_N) \in \R^{n N}$ as the control vector. Therefore, with $y = \tfrac{\partial H^{\textrm{k}}}{\partial p}$, the overall dynamics reads as:
\begin{equation*}
\left(
\begin{array}{c}
\dot{q}\\
\dot{p}
\end{array}
\right) = 
\left(
\begin{array}{cc}
0 & I_{nN}\\
0 & 0
\end{array}\right)
\left(
\begin{array}{c}
\frac{\partial H^{\textrm{k}}}{\partial q}\\
\frac{\partial H^{\textrm{k}}}{\partial p}
\end{array}
\right) + 
\left(
\begin{array}{c}
0\\
I_{nN}
\end{array}
\right) u,
\end{equation*}
where $H^{\textrm{k}}:\R^{n N} \to \R_{\geq 0}$, $H^{\textrm{k}} (q) = \sum_{i \in \mc{I}} h^{\textrm{k}}_i(\dot{q}_i) = \frac{1}{2} \dot{q}^\top \dot{q} = \frac{1}{2} p^\top p$, represents the global kinematic contribution.

\subsection{Modeling the network of agents}
By considering the pair of agents $(i, j) \in \mc{I}^2$, for instance, we say that the constraint between them is satisfied if their relative position $\|q_j(t)  - q_i(t)\|$ is lower than a critical distance $r^{\textrm{c}} > 0$, for all $t \geq 0$. 

Thus, to specify the coupling constraints among the agents within the set, we build up a graph $\mc{G} \coloneqq (V, E)$, where the vertex set $V$ coincides with $\mc{I}$ and the edge set $E$ is defined over each distance-based constraint. Consequently, $|V| = N$ and, by assuming $M$ constraints, $|E| = M$. 

Since we are dealing with a physical system and we are interested in defining a vector of relative distances, it seems reasonable to give an arbitrary orientation to each edge. Specifically, by considering the pair $(i, j) \in E$, we assume the vertex $i$ as the tail and $j$ as the head of the edge. 
Thus, the incidence matrix $B \in \R^{N \times M}$, with generic entry $b_{i,j}$, summarizes the orientation of each edge in $E$.
Therefore, the vector of relative distances $z \in \R^{n M}$ reads as
\begin{equation}\label{eq:reldistvec}
z \coloneqq (B^\top \otimes I_n) \, q = \bar{B} \, q.
\end{equation}

With a slight abuse of notation, we henceforth refer to a generic edge $j \in E$ corresponding to the pair of agents $(i, k)$. Hence, it follows that, for example, $z_{j} \coloneqq q_k - q_i$.

\subsection{Virtual couplings and energetic description of the network}
To solve Problem~1, we associate to each edge of $\mc{G}$ a virtual coupling, i.e., a pair of spring-damper in parallel. We denote with $w_j \in \R^n$ the input velocity at its ends and with $f_j \in \R^n$ the corresponding output force. Thus, the dynamics of each spring-damper subsystem reads as \cite{duindam2009modeling,van2014port}:
\begin{equation}\label{eq:virtualcoup_dynamics}
\forall j \in E:\left\{
\begin{aligned}
&\dot{z}_j = w_j,\\
&f_j = \frac{\partial h^{\textrm{s}}_{j}}{\partial z_j} + D^{\textrm{c}}_{j} w_j.
\end{aligned}
\right.
\end{equation}

Specifically, the injected damping that corresponds to each edge $j \in E$ is denoted by $D^{\textrm{c}}_{j} \in \bbS^{n}_{\succ 0}$ (dissipation matrix), while $h^{\textrm{s}}_{j}:\R^{n} \to \R_{\geq 0}$ represents the Hamiltonian contribution of the $j$-th pair of spring-damper. In particular, the latter results from a (possibly) nonlinear spring $k_j(z_j)$ and directly depends on the relative position of the agents identified by the $j$-th edge. Formally, it reads as
\begin{equation}\label{eq:elastic_potential}
h^{\textrm{s}}_{j}(z_j) = \frac{1}{2} k_j(z_j) (\|z_j\| - r_j)^2,
\end{equation}
for some rest length $0 \leq r_j \leq r^{\textrm{c}}$ and $k_j : \mc{K}_j \to \R_{> 0}$, where $\mc{K}_j \coloneqq \{z_j \in \R^n \, | \, \|z_j\| \leq  r\}$, with $r > 0$.
The choice to restrict the domain of each $k_j$ to $\mc{K}_j$, however, limits the set of possible initial positions of each agent.

Since the virtual couplings are substantially design parameters, we introduce the following assumption on the choice of the spring $k_j$, for all $j \in E$.
\smallskip
\begin{standing}
	For all $j \in E$, $k_j : \mc{K}_j \to \R_{> 0}$ is a $\mc{C}^1$, even function of the relative position $z_j$.
	{\hfill$\square$}
\end{standing}
\smallskip
\begin{figure}[!t]
	\centering
	\includegraphics[width=.9\columnwidth]{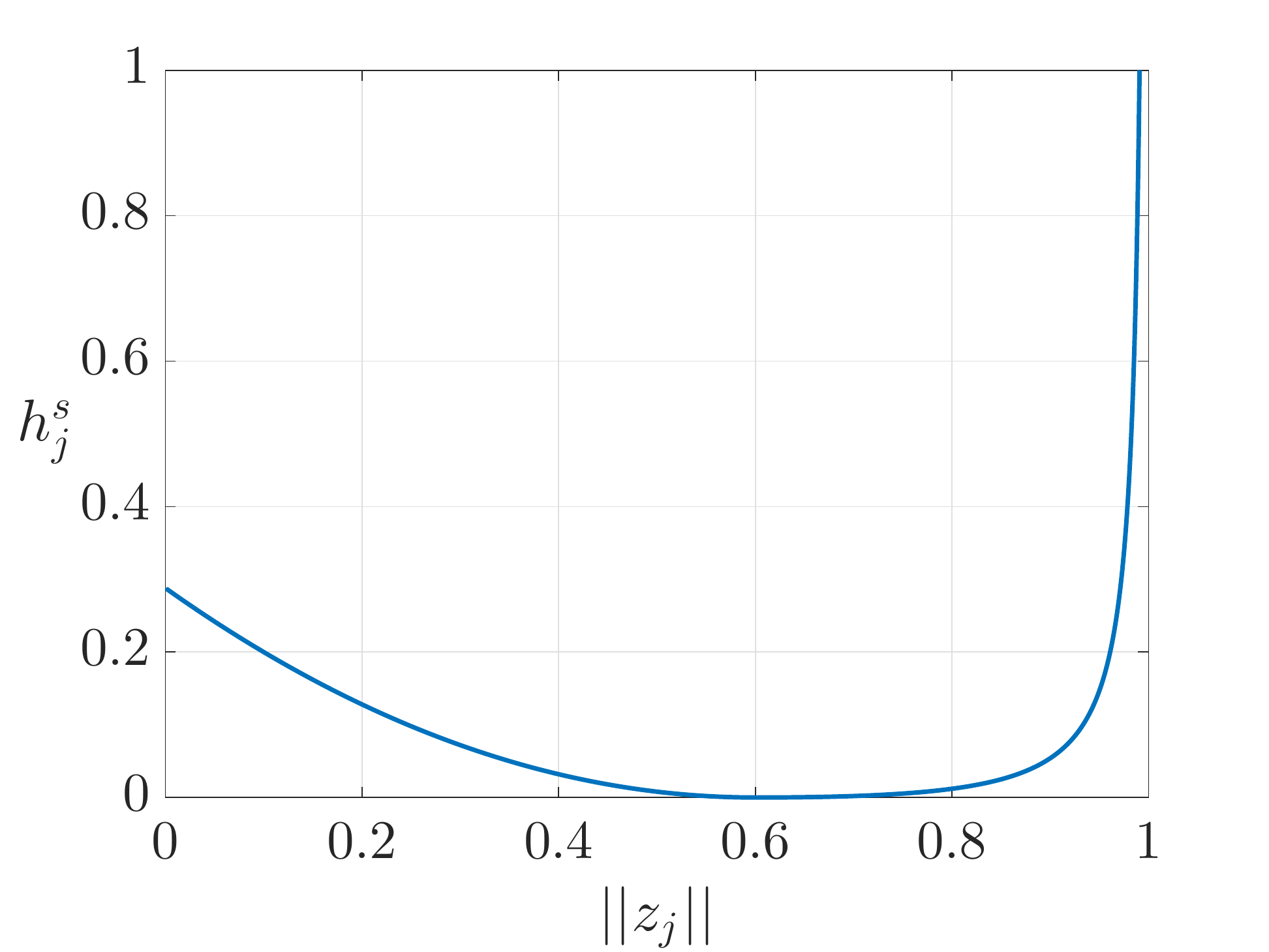}
	\caption{An example of elastic potential $h^{\textrm{s}}_j$, with $r_j = 0.6$, $r^{\textrm{c}} = 1$, $k_{j_1} = 0.8$ and $k_{j_2} = 0.06$.}
	\label{fig:ex_potcontr}
\end{figure}

In this way, $h^{\textrm{s}}_{j}$ can be shaped to obtain a desired, symmetric, intra-agents behaviour. As an example, let us consider the elastic potential shown in Fig.~\ref{fig:ex_potcontr}, which results as a contribution of a nonlinear spring, and it is formally given by:
\begin{equation}\label{eq:nonlinearspring}
h^{\textrm{s}}_{j} = \left\{
\begin{aligned}
&k_{j_1} (\|z_j\| - r_j)^2, \qquad\qquad\quad \text{if } \, \|z_j\| \leq r_j\\
&\frac{k_{j_2}}{r^{\textrm{c}} - \|z_j\|} (\|z_j\| - r_j)^2, \qquad \text{otherwise}
\end{aligned}
\right.
\end{equation}
for some $k_{j_1}$, $k_{j_2} > 0$. In this case $\mc{K}_j$ is defined by $r < r^\textrm{c}$.

Consequently, the Hamiltonian function of the whole network, $H: \R^{n M \times n N}$, is given by:
\begin{equation}\label{eq:ham_function}
\begin{aligned}
H(z, p) &= \sum_{i \in \mc{I}} h^{\textrm{k}}_i(p_i) + \sum_{j \in E} h^{\textrm{s}}_{j}(z_j)\\
& = \frac{1}{2} \sum_{i \in \mc{I}} p_i^\top p_i + \frac{1}{2} \sum_{j \in E} k_j(z_j) (\|z_j\| - r_j)^2.
\end{aligned}
\end{equation}
By choosing the vector $(z; p) \in \R^{n(M+N)}$ as state variable, we have again $y = \tfrac{\partial H}{\partial p}$ and the port-Hamilton formulation of the open-loop system is
\begin{equation}\label{eq:pH_openloop}
\left(
\begin{array}{c}
\dot{z}\\
\dot{p}
\end{array}
\right) = 
\left(
\begin{array}{cc}
0  & \bar{B}\\
0  & 0
\end{array}\right)
\left(
\begin{array}{c}
\frac{\partial H}{\partial z}\\
\vspace{-.35cm}\\
\frac{\partial H}{\partial p}
\end{array}
\right) + 
\left(
\begin{array}{c}
0\\
I_{nN}
\end{array}
\right) u.
\end{equation}

\subsection{Control synthesis and stability analysis}
The input velocity can be written as $w = \bar{B}  y = \bar{B} \tfrac{\partial H}{\partial p}$. Consequently, by exploiting the forces generated at the edges, we obtain 
the following distributed control law:
\begin{equation}\label{eq:control}
u = -\bar{B}^\top f =  - \bar{B}^\top \left( D^{\textrm{c}} \bar{B} \frac{\partial H}{\partial p} + \frac{\partial H}{\partial z}\right),
\end{equation}
with $D^{\textrm{c}} = \diag(D^{\textrm{c}}_j)$, $j \in E$.
Hence, the closed-loop system reads as:
\begin{equation}\label{eq:pH_closedloop}
\left(
\begin{array}{c}
\dot{z}\\
\dot{p}
\end{array}
\right) = 
\left(
\begin{array}{cc}
\phantom{-}0  &  \phantom{-}\bar{B}\\
-\bar{B}^\top  & -\bar{B}^\top D^{\textrm{c}} \bar{B}
\end{array}\right)
\left(
\begin{array}{c}
\frac{\partial H}{\partial z}\\
\vspace{-.35cm}\\
\frac{\partial H}{\partial p}
\end{array}
\right).
\end{equation}

Before stating a first result, we introduce the following set
\begin{equation}\label{eq:Z_set}
	\mc{Z} \coloneqq \left\{\bar{z} \in \R^{nM}  \;\Big| \; \ker(\bar{B}^\top)  \ni \frac{\partial H}{\partial z} \Big|_{z = \bar{z}}\right\},
\end{equation}
which gathers the equilibrium points among the forces generated over the edges with assigned virtual couplings.

\smallskip
\begin{proposition}\label{prop:convergence}
	The solution of the closed-loop system \eqref{eq:pH_closedloop} converges to the set of points $\mc{S} \coloneqq \{(z; p) \in \R^{n(M+N)} \, | \, z \in \mc{Z}, \, p = 0 \}$.
	{\hfill$\square$}
\end{proposition}
\smallskip
\begin{proof}
	Take $H(z, p)$ as a positive, radially unbounded, semi-definite Lyapunov candidate for the system in \eqref{eq:pH_closedloop}. Its time derivative reads as
	\begin{equation}\label{eq:ham_der}
	\dot{H} = \frac{\partial^\top H}{\partial p} \dot{p} + \frac{\partial^\top H}{\partial z} \dot{z} = - \frac{\partial^\top H}{\partial p} \bar{B}^\top D^{\textrm{c}} \bar{B} \frac{\partial H}{\partial p} \leq 0.
	\end{equation}
	Hence, by invoking the LaSalle principle, the solutions to \eqref{eq:pH_closedloop} converge to the largest invariant set where $p = 0$, i.e., 
	$$
	-\bar{B}^\top \, \frac{\partial H}{\partial z} = 0.
	$$
\end{proof}
\smallskip
\begin{corollary}\label{corollary:initcond}
	If $z_{j}(0) \in \mc{K}_j$, defined for some $r < r^\textrm{c}$, for all $j \in E$, then the solution of the closed-loop system \eqref{eq:pH_closedloop} satisfies the distance-based constraints, for all $t \geq 0$.
	{\hfill$\square$}
\end{corollary}
\smallskip
\begin{proof}
	Since $\dot{H}$ in \eqref{eq:ham_der} is negative semi-definite, given initial conditions $(z(0); p(0))$, it holds that $H(z(t), p(t)) \leq H(z(0), p(0))$ for all $t \geq 0$, and the trajectories of \eqref{eq:pH_closedloop} are bounded, namely there exists some $\mu \geq 0$, such that $\|(z(t); p(t))\| \leq \mu$ for all $t \geq 0$. 
\end{proof}
\smallskip

Note that assuming $z_{j}(0) \in \mc{K}_j$, defined for $r < r^\textrm{c}$, for all $j \in E$, is not exceedingly conservative, since distance-based constraints may be defined for communication purposes, reflecting on the possibility to compute the control in \eqref{eq:control}. Therefore, for the remainder of the paper, we consider valid such assumption.


\section{Potential game characterization}\label{sec:potgamechar}
Next, we provide some equivalence results that characterize the port-Hamiltonian system in \eqref{eq:pH_closedloop} with distance-based constraints as an exact potential game. 
For simplicity, we henceforth consider virtual springs with elastic constant independent from $z_j$, i.e., $k_j(z_j) = k_j$.

\subsection{Potential game setup}
\subsubsection{Local and collective strategies}
In our formulation, the player set coincides with $\mc{I}$, i.e., the set of double-integrator agents. Then, we assume that each agent makes decision on its position and velocity, and the local decision variable is
$$
x_i \coloneqq \left(
\begin{array}{c}
q_i\\
\dot{q}_i
\end{array} \right) \in \mc{X}_i \subset \R^{2n},
$$
for some compact and convex set $\mc{X}_i$, for all $i \in \mc{I}$. 

Clearly, since the agents follow a certain dynamics, $q_i$ and $\dot{q}_i$ can not be chosen independently. Thus, the collective vector of strategies reads as:
\begin{equation*}\label{eq:collective_vector}
\bs{x} \coloneqq
\left(
\begin{array}{c}
q\\
\dot{q}
\end{array} \right) \in \bs{\mc{X}} \subset \R^{2nN}, \; \text{where }\, \bs{\mc{X}} \coloneqq \prod_{i \in \mc{I}} \mc{X}_i.
\end{equation*}

Note that the state variable for the port-Hamiltonian system in \eqref{eq:pH_closedloop} can be written as:
\begin{equation}\label{eq:pH_var}
\bs{x_{\textrm{pH}}} \coloneqq
\left(
\begin{array}{c}
z\\
p
\end{array}
\right) =
\left(
\begin{array}{cc}
\bar{B} & 0\\
0 & I_{nN}
\end{array}
\right)
\bs{x}.
\end{equation}
\smallskip
\begin{lemma}\label{lemma:}
	The convergence of $\bs{x_{\textrm{pH}}}$ to some $\bs{\tilde{x}} \in \mc{S}$ implies the convergence of $\bs{x}$ to some $\bs{\bar{x}} \in \bs{\mc{X}}$.
	{\hfill$\square$}
\end{lemma}
\smallskip
\begin{proof}
	The point $\bs{\tilde{x}}$ is of the form $(\tilde{z}; 0)$. Then, from \eqref{eq:pH_var}, $\dot{q} \to 0$, while $q \to \bar{q}$ that satisfies $\bar{B} \bar{q} = \tilde{z}$.
\end{proof}
\smallskip

\subsubsection{Local objective functions}
By defining $\mc{L}_i \coloneqq \{j \in E \, | \, b_{i,j} \neq 0\}$ as the set of edges that involve the $i$-th player, for all $i \in \mc{I}$, we identify each local objective function as:
\begin{align}\label{eq:local_obj}
J_i(x_i,& \bs{x}_{-i}) = h_i(x_i, \bs{x}_{-i}) = h^{\textrm{k}}_i(x_i) + \sum_{\ell \in \mc{L}_i}^{} h^{\textrm{s}}_{\ell}(x_i, \bs{x}_{-i})  \nonumber\\
&= \frac{1}{2} \left(\dot{q}_i^\top \dot{q}_i + \sum_{\ell \in \mc{L}_i}^{} k_{\ell} \, (\|q_i - q_k\| - r_{\ell})^2\right).
\end{align}

Thus, we assume that every decision maker $i \in \mc{I}$ seeks for a feasible strategy that minimizes its energetic contribution within the network:
\begin{equation}\label{eq:optprob_pH}
\forall i \in \mc{I} : \underset{x_i \in \mc{X}_i}{\textrm{min}} \; J_i(x_i, \bs{x}_{-i}).
\end{equation}


\subsubsection{The potential function}
Let us consider the Hamiltonian function in \eqref{eq:ham_function}. Due to the symmetric contribution of each virtual spring, it can be equivalently rewritten as:
$$
H(z, p) = \frac{1}{2}\sum_{i \in \mc{I}}^{}\left( p_i^\top p_i + \frac{1}{2} \sum_{\ell \in \mc{L}_i}^{} k_{\ell} \, (\|z_{\ell}\| - r_{\ell})^2\right).
$$
%
\smallskip
\begin{proposition}\label{prop:equivalence}
	The game $\Gamma \coloneqq (\mc{I}, \{J_i\}_{i \in \mc{I}}, \{\mc{X}_i\}_{i \in \mc{I}})$ is an exact potential game with potential function $H(\bs{x})$.
	{\hfill$\square$}
\end{proposition}
\smallskip
\begin{proof}
	Take an arbitrary player $i \in \mc{I}$, a feasible rivals' strategy vector $\bs{x}_{-i}$ and two feasible strategies $x_i = (q_i;\dot{q}_i)$, $\bar{x}_i = (\bar{q}_i;\dot{\bar{q}}_i) \in \mc{X}_i$. Then, by directly applying \eqref{eq:exactpotgame}, we have:
	\begin{align}\label{eq:expotgame}
	&J_i(x_i, \bs{x}_{-i}) - J_i(\bar{x}_i, \bs{x}_{-i}) = \nonumber\\
	& h^\textrm{k}_i(\dot{q}_i) - h^\textrm{k}_i(\dot{\bar{q}}_i) + \sum_{\ell \in \mc{L}_i}\left( h^{\textrm{s}}_{\ell}(q_i, q_k) - h^{\textrm{s}}_{\ell}(\bar{q}_i, q_k)\right).
	\end{align}
	
	Since the kinetic contribution depends on the local variable only, by adding and subtracting in \eqref{eq:expotgame} the term $h^\textrm{k}_j(\dot{q}_j)$ for all $j \in \mc{I} \setminus \{i\}$, we obtain both $H^\textrm{k}(\dot{q})$ and $-H^\textrm{k}(\dot{\bar{q}}_i, \dot{q}_{-i})$. By referring to the elastic potential, we exploit the symmetry of the springs, separating the contribution of all players except the $i$-th one as described next:
	\begin{align}\label{eq:elpot_split}
	\frac{1}{2} \sum_{j \in \mc{I} \setminus \{i\}} \sum_{\ell \in \mc{L}_j} h^{\textrm{s}}_{\ell}(q_j,q_k) &= \frac{1}{2} \sum_{j \in \mc{I} \setminus \{i\}} \left( \sum_{\ell \in \mc{L}_j \cap \mc{L}_i} h^{\textrm{s}}_{\ell}(q_j,q_i)\right.\nonumber\\ 
	&\left.+  \sum_{\ell \in \mc{L}_j \setminus \mc{L}_i} h^{\textrm{s}}_{\ell}(q_j,q_k)\right).
	\end{align}
	
	Here, $\mc{L}_j \cap \mc{L}_i$ contains the edges that connect the agents with the $i$-th one, while $\mc{L}_j \setminus \mc{L}_i$ gathers the edges that do not directly involve the agent $i$.
	Thus, we can add and subtract \eqref{eq:elpot_split}, considering both $x_i$ and $\bar{x}_i$. In the latter case, \eqref{eq:elpot_split} is:
	$$
	\frac{1}{2} \sum_{j \in \mc{I} \setminus \{i\}} \left( \sum_{\ell \in \mc{L}_j \cap \mc{L}_i} h^{\textrm{s}}_{\ell}(q_j,\bar{q}_i) +  \sum_{\ell \in \mc{L}_j \setminus \mc{L}_i} h^{\textrm{s}}_{\ell}(q_j,q_k)\right).
	$$
	
	After some manipulations, the relation in \eqref{eq:expotgame} becomes:
	$$
	\begin{aligned}
	&J_i(x_i, \bs{x}_{-i}) - J_i(\bar{x}_i, \bs{x}_{-i}) = H^\textrm{k}(\dot{q}) -H^\textrm{k}(\dot{\bar{q}}_i, \dot{q}_{-i})\\ 
	&+ \frac{1}{2} \sum_{\ell \in \mc{L}_i} h^{\textrm{s}}_{\ell}(q_i, q_k) + \frac{1}{2} \sum_{j \in \mc{I}\setminus\{i\}} \left( \sum_{\ell \in \mc{L}_j \cap \mc{L}_i} h^{\textrm{s}}_{\ell}(q_j,q_i) \right.\\
	&\left.+ \sum_{\ell \in \mc{L}_j \setminus \mc{L}_i} h^{\textrm{s}}_{\ell}(q_j,q_k)\right) - \frac{1}{2} \sum_{\ell \in \mc{L}_i} h^{\textrm{s}}_{\ell}(\bar{q}_i, q_k) \\
	& - \frac{1}{2} \sum_{j \in \mc{I}\setminus\{i\}} \left( \sum_{\ell \in \mc{L}_j \cap \mc{L}_i} h^{\textrm{s}}_{\ell}(q_j,\bar{q}_i) + \sum_{\ell \in \mc{L}_j \setminus \mc{L}_i} h^{\textrm{s}}_{\ell}(q_j,q_k)\right)\\
	&= H(x_i, \bs{x}_{-i}) - H(\bar{x}_i, \bs{x}_{-i}),
	\end{aligned}
	$$
	which concludes the proof.
\end{proof}
\smallskip

Since each virtual coupling determines a symmetric energetic contribution, i.e., $h^{\textrm{s}}_j(q_j, q_i) = h^{\textrm{s}}_j(q_i, q_j)$, and each kinetic contribution depends on the local variable only, such a potential game belongs to the class of \gls{BSI} games \cite{monderer1996potential,ui2000shapley}. Hence, the potential function can be equivalently written as
\begin{equation}\label{eq:bsi_potfun}
H(\bs{x}) = \sum_{i \in \mc{I}} \left(h^{\textrm{k}}_i(\dot{q}_i) +  \sum_{\substack{j \in \mc{I}, \, j \prec i}} h^{\textrm{s}}_j(q_j, q_i) \right),
\end{equation}
where $j \prec i$ identifies a predefined ordering within the player set $\mc{I}$. It follows by \cite[\S2]{ui2000shapley} that every \gls{BSI} game is an exact potential game with potential function of the form \eqref{eq:bsi_potfun}. 

\subsection{Nash equilibrium seeking}
Once proved the equivalence between the control problem addressed and a potential game, here we study the convergence of the constrained, multi-agent system to some \gls{NE}.

We recall that the port-Hamiltonian model in \eqref{eq:pH_closedloop} assumes the following ``collective'', gradient-based dynamics
\begin{equation}\label{eq:collectivedynamics}
\bs{\dot{x}_{\textrm{pH}}} = 
- K \; \nabla H(\bs{x_{\textrm{pH}}}),
\end{equation}
with $ K \coloneqq \left(\begin{smallmatrix}
0  & \; & -\bar{B}\\
\bar{B}^\top & \; & \bar{B}^\top D^{\textrm{c}} \bar{B}
\end{smallmatrix}\right) \succcurlyeq 0$, and it converges to the set of points $\mc{S}$. Next, we give a convergence result, showing that each point in $\mc{S}$ corresponds to a variational equilibrium.
\smallskip
\begin{proposition}\label{prop:variationalNE}
	Any $\bs{x}^{\ast} \in \mc{S}$ is a variational equilibrium of the exact potential game $\Gamma$.
	{\hfill$\square$}
\end{proposition}
\smallskip
\begin{proof}
	For all $i \in \mc{I}$, the gradient of $J_i(x_i, \bs{x}_{-i})$ calculated with respect to the local variable $x_i$ reads as 
	$$
	\nabla_{x_i} J_i(x_i, \bs{x}_{-i}) = \left(
	\begin{array}{c}
	\underset{\ell \in \mc{L}_i}{\sum} b_{i,\ell} k_{\ell} \left(1 - \frac{r_\ell}{\|q_i - q_k\|}\right) (q_i - q_k)\\
	\vspace{-.1cm}\\
	\dot{q}_i
	\end{array} \right),
	$$
	where the element $b_{i,\ell}$ of the matrix $\bar{B}$ is needed for the correct sign of the partial derivative $\partial J_i(x_i, \bs{x}_{-i})/\partial x_i$.
	
	Thus, by stacking and re-arranging the gradient of each player, the pseudo-gradient mapping of the game reads as:
	\begin{equation*}
	\bs{F}(\bs{x}) = 
	\left(
	\begin{array}{c}
	\bar{B}^\top_{(1,:)} \frac{\partial{H}}{\partial{z}}\\
	\vdots\\
	\bar{B}^\top_{(N,:)} \frac{\partial{H}}{\partial{z}}\\
	\dot{q}_1\\
	\vdots\\
	\dot{q}_N
	\end{array}
	\right) = 
	\left(\begin{array}{cc}
	\bar{B}^\top &  0\\
	0 & I_{nN}
	\end{array}
	\right) \nabla H(\bs{x_{\textrm{pH}}}),
	\end{equation*}
	where $\bar{B}^\top_{(i,:)}$ selects the $i$-th row of the matrix $\bar{B}^\top$.
	Therefore, the pseudo-gradient $\bs{F}$ evaluated at any equilibrium point $\bs{x}^{\ast} \in \mc{S}$ is null. Since $H$ is a positive semi-definite, bounded from below function, it follows that any $\bs{x}^{\ast} \in \mc{S}$ solves the associated variational inequality problem and hence it is a variational equilibrium of the potential game.
\end{proof}
\smallskip
\smallskip
\begin{corollary}
	Given any $\bs{x}(0) \in \bs{\mc{X}}$, the closed-loop, port-Hamiltonian system in \eqref{eq:pH_closedloop} converges to a \gls{NE} of the exact potential game, satisfying the constraints for all $t \geq 0$.
	{\hfill$\square$}
\end{corollary}
\smallskip
\begin{proof}
	The proof follows directly as a consequence of Propositions~\ref{prop:convergence} and \ref{prop:variationalNE}, Lemma~\ref{lemma:} and Corollary~\ref{corollary:initcond}.
\end{proof}
\smallskip

Now, let us consider the case in which the graph $\mc{G}$ is connected and acyclic (loop-free). In this case, the control law in \eqref{eq:control} leads the system to the global minimum of $H$, i.e., at the equilibrium point in which each virtual spring is at its rest length. In view of the equivalence in Prop.~\ref{prop:equivalence}, this minimum corresponds to a \gls{NE}. Note that this condition represents the ideal outcome in several, multi-agent control problems, e.g., formation control.
\smallskip
\begin{proposition}
	Let $\mc{G}$ be a connected and acyclic graph. Then, the collective dynamics in \eqref{eq:collectivedynamics} converges to a \gls{NE} of the associated exact potential game $\Gamma$.
	{\hfill$\square$}
\end{proposition}
\smallskip
\begin{proof}
	Since the graph $\mc{G}$ is acyclic, its incidence matrix satisfies $\ker(\bar{B}^\top) = \varnothing$. Thus, by replicating the proof of Prop.~\ref{prop:convergence}, the system in \eqref{eq:pH_closedloop} converges to the largest invariant set where $p = 0$, that leads to the set of $z$ such that:
	\begin{equation*}
	\frac{\partial{H}}{\partial{z}} \Big\rvert_{z = \bar{z}} = 0.
	\end{equation*}
	This implies that $\nabla H(\bar{z}, 0) = 0$, which corresponds to the global minimum of the Hamiltonian function $H$, and hence of the exact potential function.
\end{proof}
\smallskip

Conversely, the convergence to the global minimum of $H$ does not imply that $\mc{G}$ is acyclic (see the example in \S\ref{sec:num_sim}). 

As a final remark, we stress that the convergence of $\bs{x_{\textrm{pH}}}$ in \eqref{eq:collectivedynamics} depends on $D^\textrm{c}$, containing the damping parameter of each virtual coupling, which allows to shape the transient response. In parallel, the set $\mc{Z}$ in \eqref{eq:Z_set} depends on each artificial potential $h^\textrm{s}_j$, $j \in E$, introduced with the virtual couplings. In this way, each (nonlinear) spring influences the asymptotic behaviour of the network of agents, i.e., the set of reachable equilibria of the associated potential game $\Gamma$.
\section{Numerical simulation}\label{sec:num_sim}

\begin{figure}[!h]
	\centering
	\subfloat[\label{fig:convergence}]
	{\includegraphics[width=.9\columnwidth]{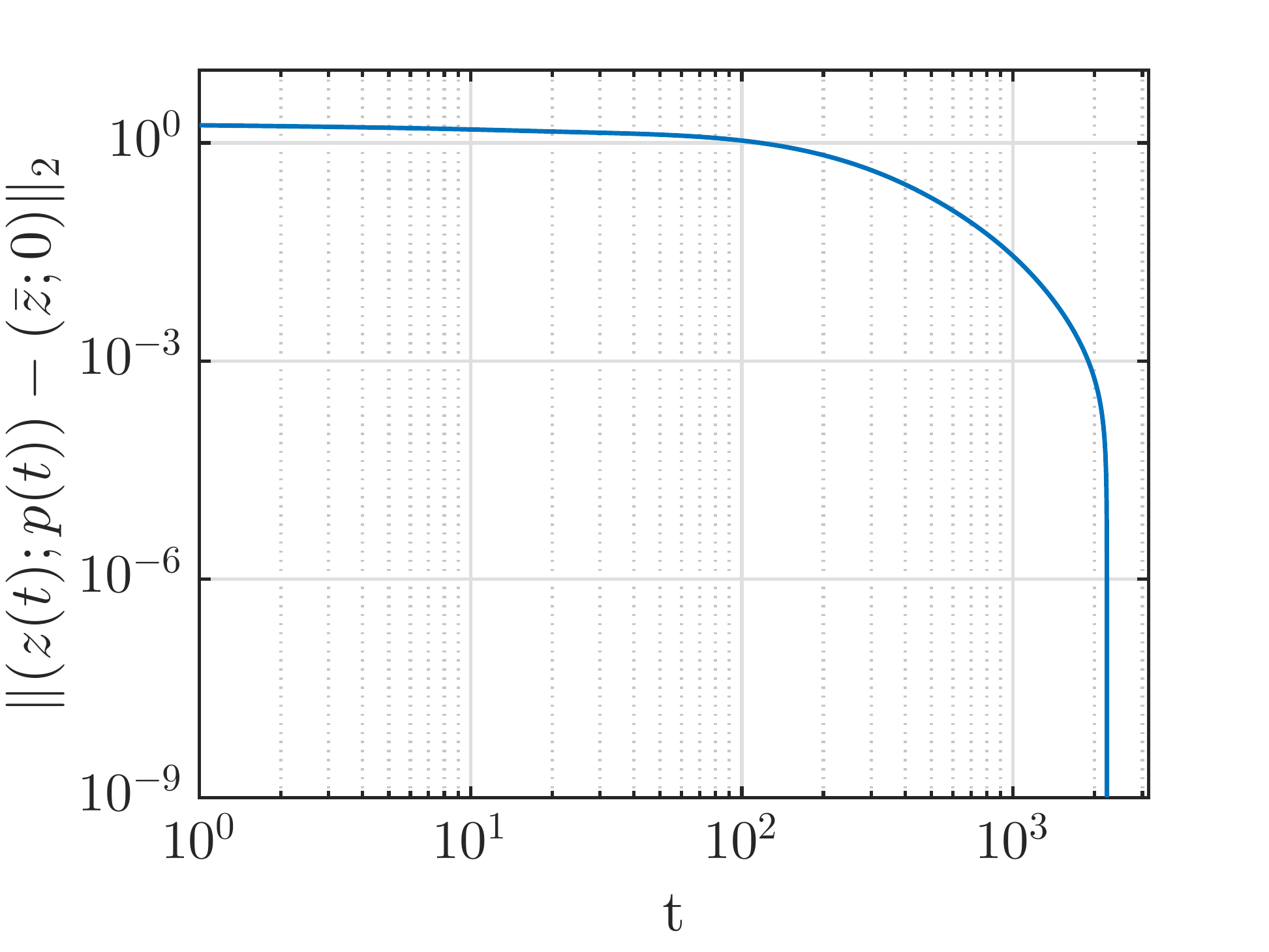}}
	\hfill
	\subfloat[\label{fig:agents_path}]
	{\includegraphics[width=.9\columnwidth]{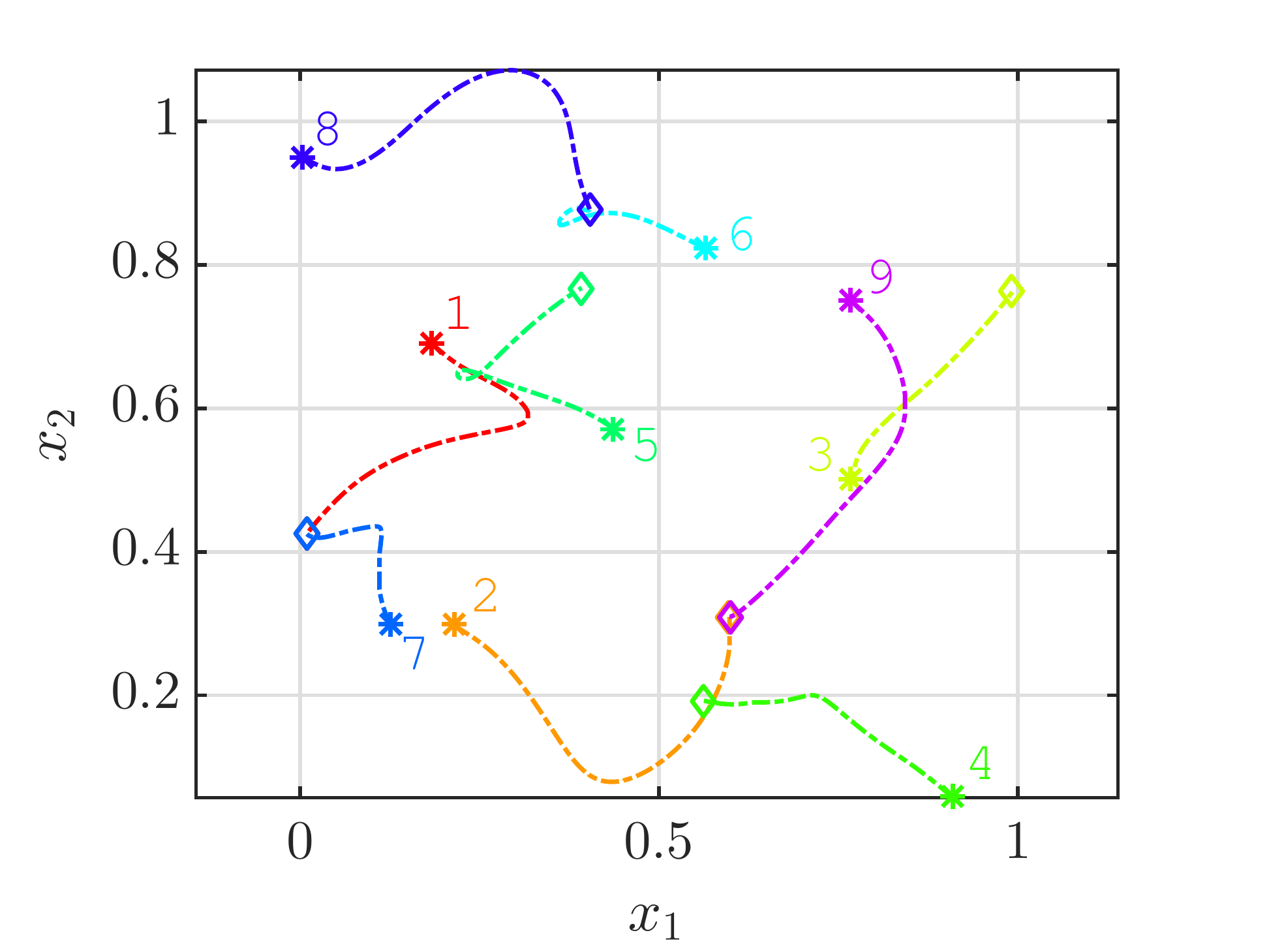}}\\
	\subfloat[\label{fig:relative_distances}]
	{\includegraphics[width=.9\columnwidth]{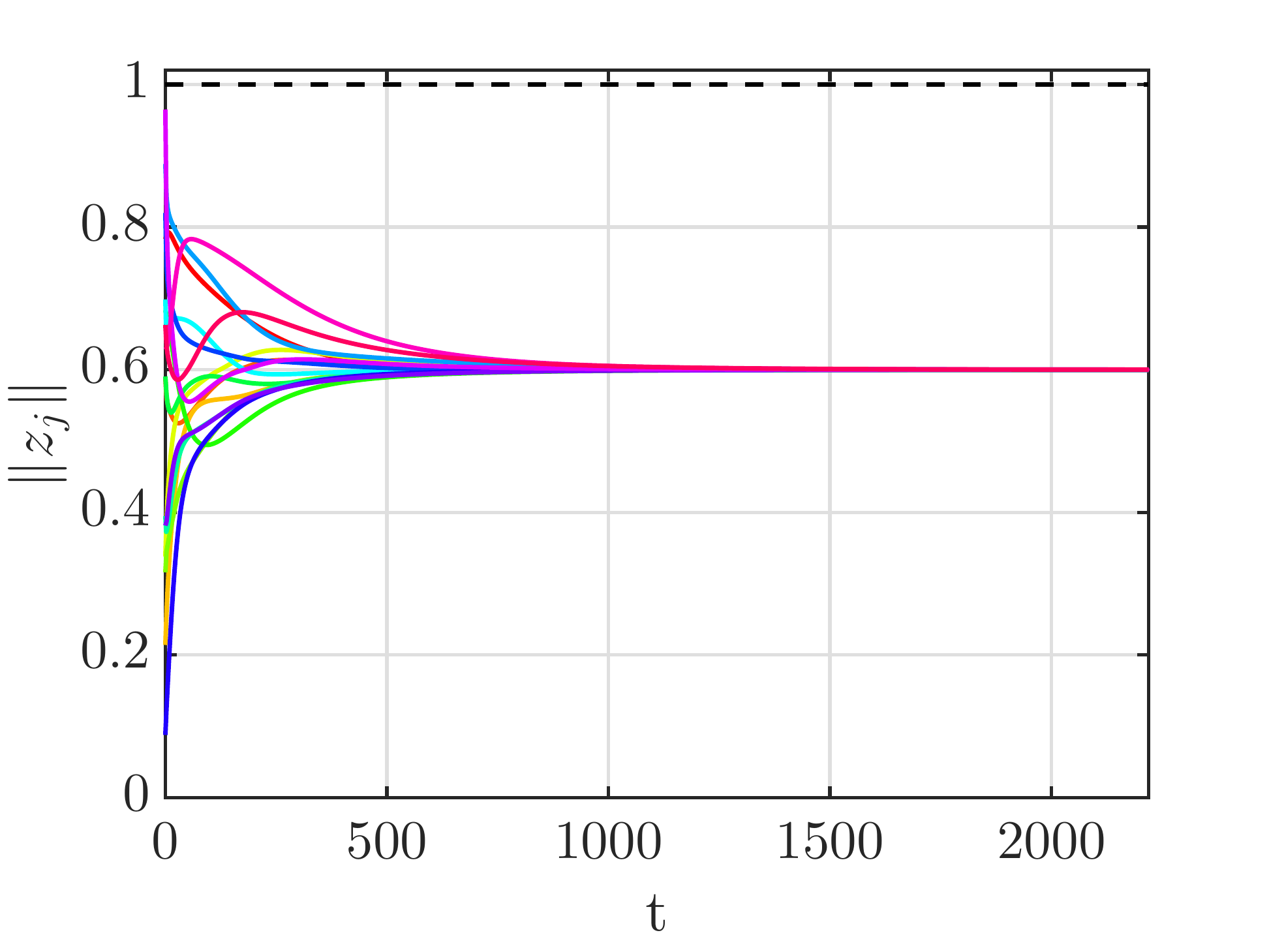}}
	\hfill
	\subfloat[\label{fig:ham_function}]
	{\includegraphics[width=.9\columnwidth]{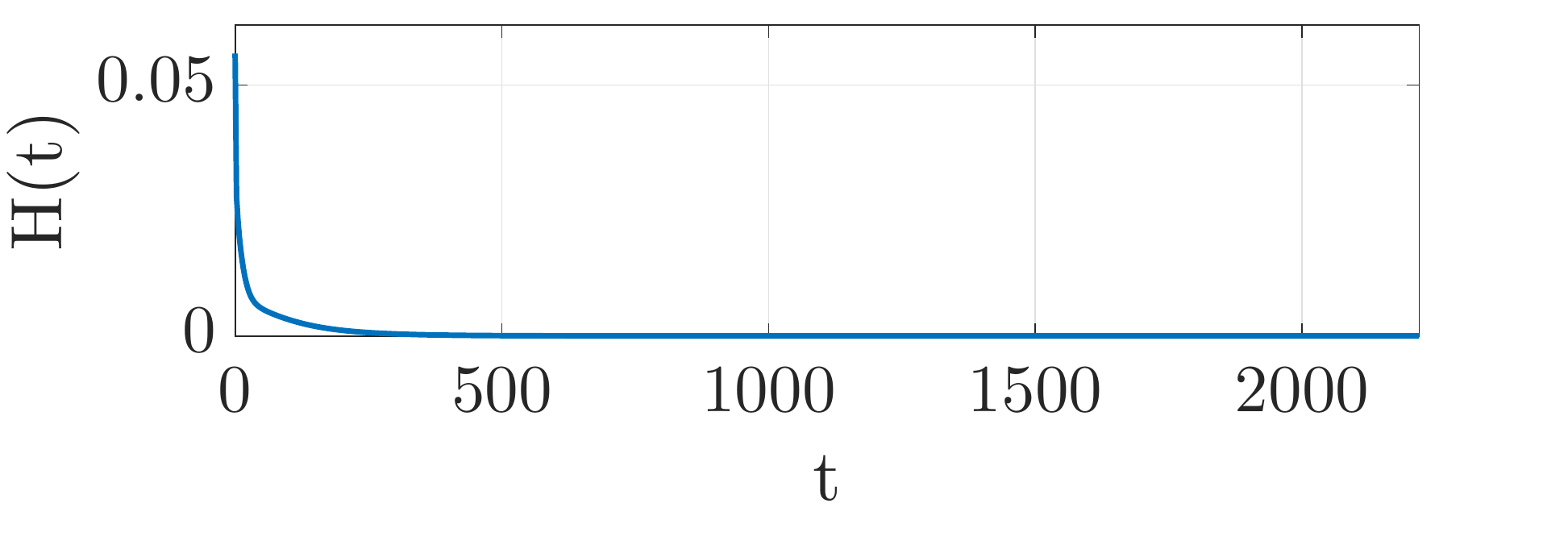}}
	\caption{(a) Convergence of $\bs{x_{\textrm{pH}}}$ to a \gls{NE}, $\bs{x}^\ast \in \mc{S}$. (b) ``Exploration'' of the constrained set of agents from a random starting point (colored asterisks) to an equilibrium (colored diamonds). (c) Relative distances. (d) Hamiltonian function.}
	\label{fig:num_sim}
\end{figure}

In Fig.~\ref{fig:num_sim} is reported an example involving $9$ agents and $16$ randomly chosen constraints. In details, we have:
\begin{equation*}
	\begin{aligned}
		E = \{&(1,2), (1,4), (1,8), (1,9), (2,3), (2,6), (2,7), (3,5),\\ 
		&(3,6), (3,8), (4,5), (4,7), (6,7), (6,9), (7,8), (7,9)\}
	\end{aligned}
\end{equation*}
A virtual coupling, designed with a constant damper and a nonlinear spring that shapes the energetic contribution in Fig.~\ref{fig:ex_potcontr}, is assigned to each edge. Fig.~\ref{fig:num_sim}(a) shows the convergence of $\bs{x_{\textrm{pH}}}$ to a \gls{NE} of the associated potential game, while Fig.~\ref{fig:num_sim}(b) highlights the ``exploration'' in $\R^2$ of the constrained set of agents, seeking for a set of positions that trade off the fulfillment of the distance-based constraints (as shown in Fig.~\ref{fig:num_sim}(c)) and the minimization of the potential function.
Note that, despite the associated graph $\mc{G}$ is cyclic, $\|z_j\| \to 0.6$ for all $j \in E$, i.e., to the rest length $r_j$ of each spring. Accordingly, $H(t) \to 0$ as $t \to \infty$ (Fig.~\ref{fig:num_sim}(d)).

\section{Conclusion and outlook}
The distributed, gradient-based control algorithm \eqref{eq:control} is suitable to steer the multi-agent system with distance-based constraints to an equilibrium point. Since the pairs of spring-damper are basically design parameters, they offer the possibility to shape the transient behaviour of the network, fulfilling the constraints and maintaining the stability. 
Moreover, by exploiting the equivalence relations with potential games, this framework allows to shape and, ideally, determine in advance the potential function, its set of local minima and, consequently, the set of \gls{NE} related with the game.

Future research will investigate an optimal procedure to shape the artificial potentials $h_j^{\textrm{s}}$, introduced with virtual couplings, and hence the reachable set of Nash equilibria. Moreover, additional technical assumptions will be investigated toward a generalization of the proposed framework, to embrace a broader class of networked, multi-agent systems that admit a port-Hamiltonian formulation.


\bibliographystyle{IEEEtran}
\bibliography{19_ECC.bib}

\end{document}